\documentclass[12pt]{article}
\usepackage{mathrsfs}
\usepackage{amssymb}
\usepackage{amsmath}
\usepackage{amsthm}
\usepackage{latexsym}
\usepackage{graphicx}
\usepackage{cases}
 \numberwithin{equation}{section}

\newcommand{\dv}{\mathrm{div}\,}
\newcommand{\cl}{\mathrm{curl}\,}

\newtheorem{Theorem}{Theorem}[section]
\newtheorem{Lemma}{Lemma}[section]

\newtheorem{Remark}{Remark}[section]
\title{\bf Blow-up criteria for Boussinesq system and MHD system and Landau-Lifshitz equations in a bounded domain}
\date{}
\author{
\bf\large Jishan Fan\\
\small Department of Applied Mathematics,\\
\small Nanjing Forestry University, Nanjing 210037, P.R.China\\
\small E-mail: fanjishan@njfu.edu.cn\\[4mm]
\bf\large Wenjun Sun,\qquad Junping Yin\thanks{Corresponding author.}\\
\small Institute of Applied Physics and Computational Mathematics,\\
\small Beijing, 100088, P.R.China\\
\small E-mail: sun\_wenjun@iapcm.ac.cn, \quad yinjp829829@126.com}
\begin{document}
\maketitle
\begin{abstract}
In this paper, we prove some blow-up criteria for the 3D Boussinesq system with zero heat conductivity and MHD system and Landau-Lifshitz equations in a bounded domain.\\[3mm]
{\bf Running Title:} Boussinesq-MHD-Landau-Lifshitz
\end{abstract}
%\newpage
\section{Introduction}

Let $\Omega$ be a bounded, simply connected domain in $\mathbb{R}^3$ with smooth boundary $\partial\Omega$, and $\nu$ be the unit outward normal vector to $\partial\Omega$. First, we consider the regularity criterion of the Boussinesq system with zero heat conductivity:
\begin{eqnarray}
&&\dv u=0,\label{1.1}\\
&&\partial_tu+u\cdot\nabla u+\nabla\pi-\Delta u=\theta e_3,\label{1.2}\\
&&\partial_t\theta+u\cdot\nabla\theta=0\ \ \mathrm{in}\ \ \Omega\times(0,\infty),\label{1.3}\\
&&u\cdot\nu=0,\cl u\times\nu=0\ \ \mathrm{on}\ \ \partial\Omega\times(0,\infty),\label{1.4}\\
&&(u,\theta)(\cdot,0)=(u_0,\theta_0)\ \ \mathrm{in}\ \ \Omega\subseteq\mathbb{R}^3,\label{1.5}
\end{eqnarray}
where $u,\pi$, and $\theta$ denote unknown velocity vector field, pressure scalar, and temperature scalar of the fluid, respectively. $\omega:=\cl u$ is the vorticity and $e_3:=(0,0,1)^t$.

When $\theta=0$, \eqref{1.1} and \eqref{1.2} are the well-known Navier-Stokes system. Giga \cite{1}, Kim \cite{2}, Kang and Kim \cite{3} have proved some Serrin type regularity criteria.

The first aim of this paper is to prove a new regularity criterion for the problem \eqref{1.1}-\eqref{1.5}, we will prove
\begin{Theorem}\label{th1.1}
Let $u_0\in H^3, \theta_0\in W^{1,p}$ with $3<p\leq 6$ and $\dv u_0=0$ in $\Omega$ and $u_0\cdot\nu=0, \cl u_0\times\nu=0$ on $\partial\Omega$. Let $(u,\theta)$ be a strong solution of the problem \eqref{1.1}-\eqref{1.5}. If $u$ satisfies
\begin{equation}
\nabla u\in L^1(0,T;BMO(\Omega))\label{1.6}
\end{equation}
with $0<T<\infty$, then the solution $(u,\theta)$ can be extended beyond $T>0$. Here $BMO$ denotes the space of bounded mean oscillation.
\end{Theorem}

Secondly, we consider the blow-up criterion of the 3D MHD system
\begin{eqnarray}
&&\dv u=\dv b=0,\label{1.7}\\
&&\partial_tu+u\cdot\nabla u+\nabla\left(\pi+\frac12|b|^2\right)-\Delta u=b\cdot\nabla b,\label{1.8}\\
&&\partial_tb+u\cdot\nabla b-b\cdot\nabla u=\Delta b\ \ in\ \ \Omega\times(0,\infty),\label{1.9}\\
&&u\cdot\nu=0,\cl u\times\nu=0,b\cdot\nu=0,\cl b\times\nu=0\ \ on\ \ \partial\Omega\times(0,\infty),\label{1.10}\\
&&(u,b)(\cdot,0)=(u_0,b_0)\ \ in\ \ \Omega\subseteq\mathbb{R}^3.\label{1.11}
\end{eqnarray}
Here $b$ is the magnetic field of the fluid.

It is well-known that the problem \eqref{1.7}-\eqref{1.11} has a unique local strong solution \cite{4}. But whether this local solution can exist globally is an outstanding problem. Kang and Kim \cite{3} prove some Serrin type regularity criteria.

The second aim of this paper is to prove a new regularity criterion for the problem \eqref{1.7}-\eqref{1.11}, we will prove
\begin{Theorem}\label{th1.2}
Let $u_0,b_0\in H^3$ with $\dv u_0=\dv b_0=0$ in $\Omega$ and $u_0\cdot\nu=b_0\cdot\nu=0, \cl u_0\times\nu=\cl b_0\times\nu=0$ on $\partial\Omega$. Let $(u,b)$ be a strong solution to the problem \eqref{1.7}-\eqref{1.11}. If \eqref{1.6} holds true, then the solution $(u,b)$ can be extended beyond $T>0$.
\end{Theorem}
\begin{Remark}\label{rel.1}
When $\Omega:=\mathbb{R}^3$, our result gives the following
well-known regularity criterion $$\omega:=\cl u\in L^1(0,T;\dot
B_{\infty,\infty}^0),$$ but the method of proof we used is different
from that in \cite{14,15}. Here $\dot B_{\infty,\infty}^0$ denotes
the homogeneous Besov space \cite{13}.
\end{Remark}

Next, we consider the following 3D density-dependent MHD equations:
\begin{eqnarray}
&&\dv u=\dv b=0,\label{1.12}\\
&&\partial_t\rho+\dv(\rho u)=0,\label{1.13}\\
&&\partial_t(\rho u)+\dv(\rho u\otimes u)+\nabla\left(\pi+\frac12|b|^2\right)-\Delta u=b\cdot\nabla b,\label{1.14}\\
&&\partial_tb+u\cdot\nabla b-b\cdot\nabla u=\Delta b\ \ in\ \ \Omega\times(0,\infty),\label{1.15}\\
&&u=0,b\cdot\nu=0,\cl b\times\nu=0\ \ on\ \ \partial\Omega\times(0,\infty),\label{1.16}\\
&&(\rho,\rho u,b)(\cdot,0)=(\rho_0,\rho_0u_0,b_0)\ \ in\ \ \Omega\subset\mathbb{R}^3.\label{1.17}
\end{eqnarray}
For this problem, in \cite{5}, Wu proved that if the initial data
$\rho_0,u_0$, and $b_0$ satisfy
\begin{equation}
0\leq\rho_0\in H^2,u_0\in H_0^1\cap H^2, b_0\in H^2, -\Delta u_0+\nabla\left(\pi_0+\frac12|b_0|^2\right)=b_0\cdot\nabla b_0+\sqrt{\rho_0}g\label{1.18}
\end{equation}
for some $(\pi_0,g)\in H^1\times L^2$, then there exists a positive time $T_*$ and a unique strong solution $(\rho, u, b)$ to the problem \eqref{1.12}-\eqref{1.17} such that
\begin{equation}
\begin{array}{l}
\rho\in C([0,T_*];H^2), u\in C([0,T_*];H_0^1\cap H^2)\cap L^2(0,T_*;H^2),\\
u_t\in L^2(0,T_*;H_0^1), \sqrt\rho u_t\in L^\infty(0,T_*;L^2),\\
b\in L^\infty(0,T_*;H^2)\cap L^2(0,T_*;H^3), b_t\in L^\infty(0,T_*;L^2)\cap L^2(0,T_*;H^1).
\end{array}\label{1.19}
\end{equation}
And when $b=0$, Kim \cite{2} proved the following regularity
criterion:
\begin{equation}
u\in L^\frac{2s}{s-3}(0,T;L_w^s(\Omega))\ \ \mathrm{with}\ \ 3<s\leq\infty.\label{1.20}
\end{equation}
Here $L_w^s$ denotes the weak-$L^s$ space and $L_w^\infty=L^\infty$.

The aim of this paper is to refine \eqref{1.20}, we will prove
\begin{Theorem}\label{th1.3}
Let $\rho_0,u_0$, and $b_0$ satisfy \eqref{1.18}. Let $(\rho, u, b)$ be a strong solution of the problem \eqref{1.12}-\eqref{1.17} in the class \eqref{1.19}. If $u$ satisfies one of the following two conditions:
\begin{eqnarray}
&&(i)\quad\int_0^T\frac{\|u(t)\|_{L_w^s}^\frac{2s}{s-3}}{1+\log(e+\|u(t)\|_{L_w^s})}dt<\infty\ \ with\ \ 3<s\leq\infty,\label{1.21}\\
&&(ii)\quad u\in L^2(0,T;BMO(\Omega))\label{1.22}
\end{eqnarray}
with $0<T<\infty$, then the solution $(\rho,u,b)$ can be extended beyond $T>0$.
\end{Theorem}

Finally, we consider the 3D Landau-Lifshitz system:
\begin{eqnarray}
&&\partial_td-\Delta d=d|\nabla d|^2+d\times\Delta d,|d|=1\ \ in\ \ \Omega\times(0,\infty),\label{1.23}\\
&&\partial_\nu d=0\ \ on\ \ \partial\Omega\times(0,\infty),\label{1.24}\\
&&d(\cdot,0)=d_0,|d_0|=1\ \ in\ \ \Omega\subseteq\mathbb{R}^3.\label{1.25}
\end{eqnarray}
Carbou and Fabrie \cite{6} showed the existence and uniqueness of local smooth solutions. When $\Omega:=\mathbb{R}^n \ (n=2,3,4)$, Fan and Ozawa \cite{7} proved some regularity criteria. The aim of this paper is to prove a logarithmic blow-up criterion for the problem \eqref{1.23}-\eqref{1.25} when $\Omega$ is a bounded domain. We will prove
\begin{Theorem}\label{th1.4}
Let $d_0\in H^3(\Omega)$ with $|d_0|=1$ in $\Omega$ and $\partial_\nu d_0=0$ on $\partial\Omega$. Let $d$ be a local smooth solution to the problem \eqref{1.23}-\eqref{1.25}. If $d$ satisfies
\begin{equation}
\int_0^T\frac{\|\nabla d\|_{L^q}^\frac{2q}{q-3}}{1+\log(e+\|\nabla d\|_{L^q})}dt<\infty\ \ with\ \ 3<q\leq\infty,\label{1.26}
\end{equation}
and $0<T<\infty$, then the solution can be extended beyond $T>0$.
\end{Theorem}

In the following section 2, we give some preliminary Lemmas which
will be used in the following sections. The proof of Theorem
\ref{th1.1} of problem \eqref{1.1} -\eqref{1.5} will be given in
section 3. The new regularly criterion of Theorem \ref{th1.2} for
the 3D MHD problem \eqref{1.7} -\eqref{1.11} will be proved in
section 4. In section 5 is the proof of the Theorem \ref{th1.3}, and
in the next section 6 we give the main proof of final Theorem
\ref{th1.4}.

%%%%%%%%%%%%%%%----------------------------------------------------------------------------------------------
\section{Preliminary Lemmas}

In the following proofs, we will use the following logarithmic
Sobolev inequality \cite{8}:
\begin{equation}
\|\nabla u\|_{L^\infty}\leq C(1+\|\nabla
u\|_{BMO}\log(e+\|u\|_{W^{s,p}}))\ \ \mathrm{with}\ \
s>1+\frac3p.\label{1.27}
\end{equation}
and the following three lemmas.
\begin{Lemma}\label{le1.1}
(\cite{9}). Let $\Omega\subseteq\mathbb{R}^3$ be a smooth bounded
domain, let $b:\Omega\rightarrow\mathbb{R}^3$ be a smooth vector
field, and let $1<p<\infty$. Then
\begin{eqnarray}
&&-\int_\Omega\Delta b\cdot b|b|^{p-2}dx=\frac12\int_\Omega|b|^{p-2}|\nabla b|^2dx+4\frac{p-2}{p^2}\int_\Omega|\nabla|b|^\frac{p}{2}|^2dx\nonumber\\
&&\qquad-\int_{\partial\Omega}|b|^{p-2}(b\cdot\nabla)b\cdot\nu
d\sigma-\int_{\partial\Omega}|b|^{p-2}(\cl b\times\nu)\cdot b
d\sigma.\label{1.28}
\end{eqnarray}
\end{Lemma}
\begin{Lemma}\label{le1.2}
(\cite{10,11}). Let $\Omega$ be a smooth and bounded open set and
let $1<p<\infty$. Then the following estimate:
\begin{equation}
\|b\|_{L^p(\partial\Omega)}\leq
C\|b\|_{L^p(\Omega)}^{1-\frac1p}\|b\|_{W^{1,p}(\Omega)}^\frac{1}{p}\label{1.29}
\end{equation}
holds for any $b\in W^{1,p}(\Omega)$.
\end{Lemma}
\begin{Lemma}\label{le1.3}
There holds
\begin{equation}
\|f\|_{L^\infty(\Omega)}\leq
C(1+\|f\|_{BMO(\Omega)}\log^\frac12(e+\|f\|_{W^{1,4}(\Omega)}))\label{1.30}
\end{equation}
for any $f\in W_0^{1,4}(\Omega)$.
\end{Lemma}

\noindent{\bf Proof.} { When $\Omega:=\mathbb{R}^3$, \eqref{1.30}
has been proved in Ogawa \cite{12}. When $\Omega$ is a bounded
domain in $\mathbb{R}^3$. We can define $$\tilde
f:=\left\{\begin{array}{lcl}
f&\quad in\quad&\Omega,\\
0&\quad in\quad&\Omega^c:=\mathbb{R}^3\setminus\Omega.
\end{array}\right.$$

Then we have [10, Page 71] $$\|\tilde
f\|_{W^{1,4}(\mathbb{R}^3)}=\|f\|_{W^{1,4}(\Omega)}$$ and it is
obvious that $$\|\tilde
f\|_{L^\infty(\mathbb{R}^3)}=\|f\|_{L^\infty(\Omega)}, \|\tilde
f\|_{BMO(\mathbb{R}^3)}=\|f\|_{BMO(\Omega)}.$$

Thus \eqref{1.30} is proved.}

\hfill$\square$

Finally, when $b$ satisfies $b\cdot\nu=0$ on $\partial\Omega$, we
will also use the identity
\begin{equation}
(b\cdot\nabla)b\cdot\nu=-(b\cdot\nabla)\nu\cdot b\ \ \mathrm{on}\ \
\partial\Omega\label{1.31}
\end{equation}
for any sufficiently smooth vector field $b$.
\section{Proof of Theorem \ref{th1.1}}

This section is devoted to the proof of Theorem \ref{th1.1}. Since it is easy to prove that the problem \eqref{1.1}
-\eqref{1.5} has a unique local-in-time strong solution, we omit the details here. We only need to establish a priori estimates.

First, thanks to the maximum principle, it follows from \eqref{1.1} and \eqref{1.3} that
\begin{equation}
\|\theta\|_{L^\infty(0,T;L^\infty)}\leq C.\label{2.1}
\end{equation}

Testing \eqref{1.2} by $u$ and using \eqref{1.1} and \eqref{2.1}, we see that $$\frac12\frac{d}{dt}\int_\Omega u^2 dx+\int_\Omega|\cl u|^2dx\leq\int_\Omega\theta e_3\cdot u dx\leq\frac12\int_\Omega\theta^2dx+\frac12\int_\Omega u^2 dx,$$ which gives
\begin{equation}
\|u\|_{L^\infty(0,T;L^2)}+\|u\|_{L^2(0,T;H^1)}\leq C.\label{2.2}
\end{equation}

Applying $\cl$ to \eqref{1.2} and setting $\omega:=\cl u$, we find that
\begin{equation}
\partial_t\omega+u\cdot\nabla\omega-\Delta\omega=\omega\cdot\nabla u+\cl(\theta e_3).\label{2.3}
\end{equation}

Testing \eqref{2.3} by $\omega$ and using \eqref{1.1} and \eqref{2.1}, we infer that
\begin{eqnarray*}
\frac12\frac{d}{dt}\int_\Omega|\omega|^2dx+\int_\Omega|\cl\omega|^2dx&=&\int_\Omega(\omega\cdot\nabla)u\cdot\omega dx+\int_\Omega\theta e_3\cl\omega dx\\
&\leq&\|\nabla u\|_{L^\infty}\int_\Omega\omega^2dx+\frac12\int_\Omega|\cl\omega|^2dx+C,
\end{eqnarray*}
which implies
\begin{eqnarray*}
\frac{d}{dt}\int_\Omega|\omega|^2dx+\int_\Omega|\cl\omega|^2dx&\leq&C\|\nabla u\|_{L^\infty}\int_\Omega|\omega|^2dx+C\\
&\leq&C(1+\|\nabla u\|_{BMO})\log(e+\|u\|_{H^3})\int_\Omega|\omega|^2dx+C,
\end{eqnarray*}
and therefore
\begin{equation}
\int_\Omega|\omega|^2dx+\int_{t_0}^t\|\cl\omega\|_{L^2}^2d\tau\leq C(e+y)^{C_0\epsilon}\label{2.4}
\end{equation}
provided that
\begin{equation}
\int_{t_0}^t\|\nabla u\|_{BMO}d\tau\leq\epsilon<<1\label{2.5}
\end{equation}
and $y(t):=\sup\limits_{[t_0,t]}\|u\|_{H^3}$ for any $0<t_0\leq t\leq T$ and $C_0$ is an absolute constant.

Applying $\partial_t$ to \eqref{1.2}, we deduce that
\begin{equation}
\partial_t^2u+u\cdot\nabla u_t+\nabla\pi_t-\Delta u_t=-u_t\cdot\nabla u+\theta_t e_3.\label{2.6}
\end{equation}

Testing \eqref{2.6} by $u_t$, using \eqref{1.1}, \eqref{1.3}, \eqref{2.1} and \eqref{2.2}, we derive
\begin{eqnarray*}
&&\frac12\frac{d}{dt}\int_\Omega|u_t|^2dx+\int_\Omega|\cl u_t|^2dx\\
&=&-\int_\Omega u_t\cdot\nabla u\cdot u_t dx+\int_\Omega\theta_te_3u_t dx\\
&=&-\int_\Omega u_t\cdot\nabla u\cdot u_t dx-\int_\Omega\dv(u\theta)e_3u_t dx\\
&=&-\int_\Omega u_t\cdot\nabla u\cdot u_t dx+\int_\Omega u\theta\nabla(e_3u_t)dx\\
&\leq&\|\nabla u\|_{L^\infty}\int_\Omega|u_t|^2dx+\frac12\int_\Omega|\cl u_t|^2dx+C\\
&\leq&C(1+\|\nabla u\|_{BMO})\log(e+y)\int_\Omega|u_t|^2dx+\frac12\int_\Omega|\cl u_t|^2dx+C,
\end{eqnarray*}
which yields
\begin{equation}
\int_\Omega|u_t|^2dx+\int_{t_0}^t\int_\Omega|\cl u_t|^2dx d\tau\leq C(e+y)^{C_0\epsilon}.\label{2.7}
\end{equation}

On the other hand, thanks to the $H^2$-theory of the Stokes system, if follows from \eqref{1.2}, \eqref{2.1}, \eqref{2.4} and \eqref{2.7} that
\begin{eqnarray*}
\|u\|_{H^2}&\leq&C\|-\Delta u+\nabla\pi\|_{L^2}\\
&\leq&C\|\partial_tu+u\cdot\nabla u-\theta e_3\|_{L^2}\\
&\leq&C\|u_t\|_{L^2}+C\|u\|_{L^6}\|\nabla u\|_{L^3}+C\|\theta\|_{L^2}\\
&\leq&C\|u_t\|_{L^2}+C\|\nabla u\|_{L^2}^\frac32\|u\|_{H^2}^\frac12+C,
\end{eqnarray*}
which implies
\begin{equation}
\|u\|_{H^2}\leq C\|u_t\|_{L^2}+C\|\nabla u\|_{L^2}^3+C\leq C(e+y)^{C_0\epsilon}.\label{2.8}
\end{equation}

Applying $\nabla$ to \eqref{1.3}, testing by $|\nabla\theta|^{p-2}\nabla\theta\ (2\leq p<\infty)$ and using \eqref{1.1}, we get
\begin{eqnarray*}
\frac{d}{dt}\|\nabla\theta\|_{L^p}&\leq&C\|\nabla u\|_{L^\infty}\|\nabla\theta\|_{L^p}\\
&\leq&C(1+\|\nabla u\|_{BMO})\log(e+y)\|\nabla\theta\|_{L^p},
\end{eqnarray*}
which leads to
\begin{equation}
\|\nabla\theta\|_{L^\infty(t_0,t;L^p)}\leq C(e+y)^{C_0\epsilon}\ \ \mathrm{with}\ \ 2\leq p<\infty.\label{2.9}
\end{equation}

Testing \eqref{2.6} by $-\Delta u_t+\nabla\pi_t$, using \eqref{1.1}, \eqref{1.3}, \eqref{2.7}, \eqref{2.8} and \eqref{2.9}, we obtain
\begin{eqnarray*}
&&\frac12\frac{d}{dt}\int_\Omega|\cl u_t|^2dx+\int_\Omega|-\Delta u_t+\nabla\pi_t|^2dx\\
&=&\int_\Omega(-u_t\cdot\nabla u+\theta_te_3-u\cdot\nabla u_t)(-\Delta u_t+\nabla\pi_t)dx\\
&\leq&(\|\nabla u\|_{L^6}\|u_t\|_{L^3}+\|u\|_{L^\infty}\|\nabla\theta\|_{L^2}+\|u\|_{L^\infty}\|\nabla u_t\|_{L^2})\|-\Delta u_t+\nabla\pi_t\|_{L^2}\\
&\leq&\|u\|_{H^2}(\|u_t\|_{H^1}+\|\nabla\theta\|_{L^2})\|-\Delta u_t+\nabla\pi_t\|_{L^2}\\
&\leq&\frac12\|-\Delta u_t+\nabla\pi_t\|_{L^2}^2+C\|u\|_{H^2}^2(\|u_t\|_{H^1}^2+\|\nabla\theta\|_{L^2}^2),
\end{eqnarray*}
which leads to
\begin{equation}
\int_\Omega|\cl u_t|^2dx+\int_{t_0}^t\|u_t\|_{H^2}^2d\tau\leq C(e+y)^{C_0\epsilon}.\label{2.10}
\end{equation}

On the other hand, if follows from \eqref{2.3}, \eqref{2.10}, \eqref{2.9} and \eqref{2.8} that
\begin{eqnarray*}
\|u\|_{H^3}&\leq&C(1+\|\Delta\omega\|_{L^2})\\
&\leq&C(1+\|\partial_t\omega+u\cdot\nabla\omega-\omega\cdot\nabla u-\cl(\theta e_3)\|_{L^2})\\
&\leq&C(1+\|\partial_t\omega\|_{L^2}+\|u\|_{L^\infty}\|\nabla\omega\|_{L^2}+\|\omega\|_{L^4}\|\nabla u\|_{L^4}+\|\nabla\theta\|_{L^2})\\
&\leq&C(e+y)^{C_0\epsilon},
\end{eqnarray*}
which gives
\begin{equation}
\|u\|_{L^\infty(0,T;H^3)}\leq C,\label{2.11}
\end{equation}
and
\begin{equation}
\|\theta\|_{L^\infty(0,T;W^{1,p})}\leq C\ \ \mathrm{with}\ \ 3\leq p\leq 6.\label{2.12}
\end{equation}

This completes the proof of Theorem \ref{th1.1}.

\hfill$\square$

%%%%%%%%%%%%%%%----------------------------------------------------------------------------------------------
\section{Proof of Theorem \ref{th1.2}}

This section is devoted to the proof of Theorem \ref{th1.2}, we only need to prove a priori estimates.

First, testing \eqref{1.8} by $u$ and using \eqref{1.7}, we see that
\begin{equation}
\frac12\frac{d}{dt}\int_\Omega u^2dx+\int_\Omega|\cl u|^2dx=\int_\Omega(b\cdot\nabla)b\cdot u dx.\label{3.1}
\end{equation}

Testing \eqref{1.9} by $b$ and using \eqref{1.7}, we find that
\begin{equation}
\frac12\frac{d}{dt}\int_\Omega b^2 dx+\int_\Omega|\cl b|^2dx=\int_\Omega(b\cdot\nabla)u\cdot b dx.\label{3.2}
\end{equation}

Summing up \eqref{3.1} and \eqref{3.2}, we get the well-known energy inequality
\begin{equation}
\frac12\frac{d}{dt}\int_\Omega(u^2+b^2)dx+\int_\Omega(|\cl u|^2+|\cl b|^2)dx\leq 0.\label{3.3}
\end{equation}

Testing \eqref{1.9} by $|b|^{p-2}b\ (2\leq p\leq 6)$, using \eqref{1.7}, \eqref{1.28}, \eqref{1.29} and \eqref{1.31}, we derive
\begin{eqnarray*}
&&\frac1p\frac{d}{dt}\int_\Omega|b|^p dx+\frac12\int_\Omega|b|^{p-2}|\nabla b|^2dx+4\frac{p-2}{p^2}\int_\Omega|\nabla|b|^\frac{p}{2}|^2dx\\
&=&-\int_{\partial\Omega}|b|^{p-2}(b\cdot\nabla)\nu\cdot b d\sigma+\int_\Omega b\cdot\nabla u\cdot|b|^{p-2}b dx\\
&\leq&C\int_{\partial\Omega}|b|^p dx+\|\nabla u\|_{L^\infty}\int_\Omega|b|^p dx\\
&\leq&2\frac{p-2}{p^2}\int_\Omega|\nabla|b|^\frac{p}{2}|^2dx+C(1+\|\nabla u\|_{L^\infty})\int_\Omega|b|^p dx\\
&\leq&2\frac{p-2}{p^2}\int_\Omega|\nabla|b|^\frac{p}{2}|^2dx+C(1+\|\nabla u\|_{BMO})\int_\Omega|b|^p dx\log(e+y),
\end{eqnarray*}
which implies
\begin{equation}
\|b\|_{L^\infty(t_0,t;L^p)}+\int_{t_0}^t\int_\Omega|b|^2|\nabla b|^2dx d\tau\leq C(e+y)^{C_0\epsilon}\ \ \mathrm{with}\ \ 2\leq p\leq 6,\label{3.4}
\end{equation}
with the same $y$ and $\epsilon$ as that in \eqref{2.5}.

Taking $\cl$ to \eqref{1.8} and \eqref{1.9}, respectively, and setting $\omega:=\cl u$ and $j:=\cl b$, we infer that
\begin{eqnarray}
&&\partial_t\omega+u\cdot\nabla\omega-\Delta\omega=\omega\cdot\nabla u+b\cdot\nabla j+\sum\limits_i\nabla b_i\times\partial_ib,\label{3.5}\\
&&\partial_tj+u\cdot\nabla j-\Delta j=b\cdot\nabla\omega+\sum\limits_i\nabla b_i\times\partial_iu-\sum\limits_i\nabla u_i\times\partial_ib.\label{3.6}
\end{eqnarray}

Testing \eqref{3.5} and \eqref{3.6} by $\omega$ and $j$, respectively, summing up the result and using \eqref{1.7}, we have
\begin{eqnarray*}
&&\frac12\frac{d}{dt}\int_\Omega(\omega^2+j^2)dx+\int_\Omega(|\cl\omega|^2+|\cl j|^2)dx\\
&=&\int_\Omega(\omega\cdot\nabla)u\cdot\omega dx+\sum\limits_i\int_\Omega(\nabla b_i\times\partial_ib)\omega dx+\sum\limits_i\int_\Omega(\nabla b_i\times\partial_iu)\cdot j dx-\sum\limits_i\int_\Omega(\nabla u_i\times\partial_ib)\cdot j dx\\
&\leq&C\|\nabla u\|_{L^\infty}\int_\Omega(\omega^2+j^2)dx\\
&\leq&C(1+\|\nabla u\|_{BMO})\int_\Omega(\omega^2+j^2)dx\log(e+y),
\end{eqnarray*}
which implies
\begin{equation}
\int_\Omega(\omega^2+j^2)dx+\int_{t_0}^t\int_\Omega(|\cl\omega|^2+|\cl j|^2)dx d\tau\leq C(e+y)^{C_0\epsilon}.\label{3.7}
\end{equation}
Thus, it follows from \eqref{1.8}, \eqref{1.9} and \eqref{3.7} that
\begin{equation}
\int_{t_0}^t\int_\Omega(|u_t|^2+|b_t|^2)dx d\tau\leq C(e+y)^{C_0\epsilon}.\label{3.8}
\end{equation}

Applying $\partial_t$ to \eqref{1.8}, we have
\begin{equation}
\partial_t^2u+u\cdot\nabla u_t+\nabla\pi_t-\Delta u_t=\dv(b\otimes b)_t-u_t\cdot\nabla u.\label{3.9}
\end{equation}
Testing \eqref{3.9} by $u_t$ and using \eqref{1.7}, we get
\begin{eqnarray}
&&\frac12\frac{d}{dt}\int_\Omega|u_t|^2dx+\int_\Omega|\cl u_t|^2dx\nonumber\\
&=&-\sum\limits_{i,j}\int_\Omega(b^ib^j)_t\partial_ju_t^idx-\int_\Omega u_t\cdot\nabla u\cdot u_t dx\nonumber\\
&\leq&C\|b_t\|_{L^3}\|b\|_{L^6}\|\nabla u_t\|_{L^2}+\|\nabla u\|_{L^2}\|u_t\|_{L^4}^2\nonumber\\
&\leq&C\|b_t\|_{L^2}^\frac12\|\cl b_t\|_{L^2}^\frac12\|\cl u_t\|_{L^2}\|b\|_{L^6}+C\|\nabla u\|_{L^2}\|u_t\|_{L^2}^\frac12\|\cl u_t\|_{L^2}^\frac32\nonumber\\
&\leq&\delta\|\cl u_t\|_{L^2}^2+\delta\|\cl b_t\|_{L^2}^2+C\|b_t\|_{L^2}^2\|b\|_{L^6}^4+C\|\nabla u\|_{L^2}^4\|u_t\|_{L^2}^2\label{3.10}
\end{eqnarray}
for any $\delta\in(0,1)$.

Applying $\partial_t$ to \eqref{1.9}, we have
\begin{equation}
\partial_t^2b+u\cdot\nabla b_t-\Delta b_t=b_t\cdot\nabla u+b\cdot\nabla u_t-u_t\cdot\nabla b.\label{3.11}
\end{equation}
Testing \eqref{3.11} by $b_t$ and using \eqref{1.7}, we deduce that
\begin{eqnarray}
&&\frac12\frac{d}{dt}\int_\Omega|b_t|^2dx+\int_\Omega|\cl b_t|^2dx\nonumber\\
&=&\int_\Omega(b_t\cdot\nabla u+b\cdot\nabla u_t-u_t\cdot\nabla b)b_t dx\nonumber\\
&\leq&\|\nabla u\|_{L^2}\|b_t\|_{L^4}^2+\|b\|_{L^6}\|\nabla u_t\|_{L^2}\|b_t\|_{L^3}+\|\nabla b\|_{L^2}\|u_t\|_{L^4}\|b_t\|_{L^4}\nonumber\\
&\leq&\delta\|\cl b_t\|_{L^2}^2+\delta\|\cl u_t\|_{L^2}^2\nonumber\\
&&+C\|\nabla u\|_{L^2}^4\|b_t\|_{L^2}^2+C\|b\|_{L^6}^4\|b_t\|_{L^2}^2+C\|\nabla b\|_{L^2}^4(\|u_t\|_{L^2}^2+\|b_t\|_{L^2}^2)\label{3.12}
\end{eqnarray}
for any $\delta\in(0,1)$.

Combining \eqref{3.10} and \eqref{3.12} and taking $\delta$ small enough and using \eqref{3.7} and \eqref{3.8}, we have
\begin{equation}
\int_\Omega(|u_t|^2+|b_t|^2)dx+\int_{t_0}^t\int_\Omega(|\cl u_t|^2+|\cl b_t|^2)dx d\tau\leq C(e+y)^{C_0\epsilon}.\label{3.13}
\end{equation}
It follows from \eqref{1.8}, \eqref{1.9}, \eqref{3.7} and \eqref{3.13} that
\begin{equation}
\|u\|_{L^\infty(t_0,t;H^2)}+\|b\|_{L^\infty(t_0,t;H^2)}\leq C(e+y)^{C_0\epsilon}.\label{3.14}
\end{equation}

Testing \eqref{3.9} by $\nabla\left(\pi+\frac12|b|^2\right)_t-\Delta u_t$, and using \eqref{1.7}, we find that
\begin{eqnarray}
&&\frac12\frac{d}{dt}\int_\Omega|\cl u_t|^2dx+\int_\Omega\left|\nabla\left(\pi+\frac12|b|^2\right)_t-\Delta u_t\right|^2dx\nonumber\\
&=&\int_\Omega((b\cdot\nabla b)_t-u_t\cdot\nabla u-u\cdot\nabla u_t)\left(\nabla\left(\pi+\frac12|b|^2\right)_t-\Delta u_t\right)dx\nonumber\\
&\leq&C(\|b\|_{L^\infty}\|\nabla b_t\|_{L^2}+\|b_t\|_{L^6}\|\nabla b\|_{L^3}+\|u_t\|_{L^6}\|\nabla u\|_{L^3}\nonumber\\
&&+\|u\|_{L^\infty}\|\nabla u_t\|_{L^2})\left\|\nabla\left(\pi+\frac12|b|^2\right)_t-\Delta u_t\right\|_{L^2}\nonumber\\
&\leq&\frac14\left\|\nabla\left(\pi+\frac12|b|^2\right)_t-\Delta u_t\right\|_{L^2}^2+C(\|u\|_{L^\infty}^2+\|\nabla u\|_{L^3}^2)\|\nabla u_t\|_{L^2}^2\nonumber\\
&&+C(\|b\|_{L^\infty}^2+\|\nabla b\|_{L^3}^2)\|\nabla b_t\|_{L^2}^2.\label{3.15}
\end{eqnarray}
Similarly, testing \eqref{3.11} by $-\Delta b_t$, we infer that
\begin{eqnarray}
&&\frac12\frac{d}{dt}\int_\Omega|\cl b_t|^2dx+\int_\Omega|\Delta b_t|^2dx\nonumber\\
&=&\int_\Omega(u_t\cdot\nabla b+u\cdot\nabla b_t-b_t\cdot\nabla u-b\cdot\nabla u_t)\Delta b_t dx\nonumber\\
&\leq&(\|u_t\|_{L^6}\|\nabla b\|_{L^3}+\|u\|_{L^\infty}\|\nabla b_t\|_{L^2}+\|\nabla u\|_{L^3}\|b_t\|_{L^6}+\|b\|_{L^\infty}\|\nabla u_t\|_{L^2})\|\Delta b_t\|_{L^2}\nonumber\\
&\leq&\frac14\|\Delta b_t\|_{L^2}^2+C(\|u\|_{L^\infty}^2+\|\nabla u\|_{L^3}^2)\|\nabla b_t\|_{L^2}^2+C(\|b\|_{L^\infty}^2+\|\nabla b\|_{L^3}^2)\|\nabla u_t\|_{L^2}^2.\label{3.16}
\end{eqnarray}
Combining \eqref{3.15} and \eqref{3.16} and using \eqref{3.14} and \eqref{3.13}, we have
\begin{equation}
\int_\Omega(|\cl u_t|^2+|\cl b_t|^2)dx+\int_{t_0}^t\int_\Omega(|\Delta u_t|^2+|\Delta b_t|^2)dx d\tau\leq C(e+y)^{C_0\epsilon}.\label{3.17}
\end{equation}

On the other hand, it follows from \eqref{3.5}, \eqref{3.6}, \eqref{3.3}, \eqref{3.17} and \eqref{3.14} that
\begin{eqnarray*}
&&\|u(t)\|_{H^3}+\|b(t)\|_{H^3}\leq C(1+\|\Delta\omega\|_{L^2}+\|\Delta j\|_{L^2})\\
&\leq&C(1+\|\partial_t\omega+u\cdot\nabla\omega-\omega\cdot\nabla u-b\cdot\nabla j-\sum\limits_i\nabla b_i\times\partial_ib\|_{L^2}\\
&&+\|\partial_tj+u\cdot\nabla j-b\cdot\nabla\omega+\sum\limits_i\nabla u_i\times\partial_ib-\sum\limits_i\nabla b_i\times\partial_iu\|_{L^2})\\
&\leq&C(e+y(t))^{C_0\epsilon},
\end{eqnarray*}
which yields $$\|u\|_{L^\infty(0,T;H^3)}+\|b\|_{L^\infty(0,T;H^3)}\leq C,$$

This completes the proof of Theorem \ref{th1.2}.

\hfill$\square$

%%%%%%%%%%%%%%%----------------------------------------------------------------------------------------------
\section{Proof of Theorem \ref{th1.3}}
This section is devoted to the proof of Theorem \ref{th1.3}, we only need to establish a priori estimates.

First, it follows from \eqref{1.12} and \eqref{1.13} that
\begin{equation}
\|\rho\|_{L^\infty(0,T;L^\infty)}\leq C.\label{4.1}
\end{equation}
Testing \eqref{1.14} by $u$ and using \eqref{1.12} and \eqref{1.13}, we see that
\begin{equation}
\frac12\frac{d}{dt}\int_\Omega\rho u^2dx+\int_\Omega|\nabla u|^2dx=\int_\Omega(b\cdot\nabla)b\cdot u dx.\label{4.2}
\end{equation}
And testing \eqref{1.15} by $b$ and using \eqref{1.12} and
\eqref{1.16}, we find that
\begin{equation}
\frac12\frac{d}{dt}\int_\Omega|b|^2dx+\int_\Omega|\cl b|^2dx=\int_\Omega(b\cdot\nabla)u\cdot b dx.\label{4.3}
\end{equation}
Summing up \eqref{4.2} and \eqref{4.3}, we get the well-known energy inequality
\begin{equation}
\frac12\frac{d}{dt}\int_\Omega(\rho|u|^2+|b|^2)dx+\int_\Omega(|\nabla u|^2+|\cl b|^2)dx\leq 0.\label{4.4}
\end{equation}

\noindent(I) Let \eqref{1.21} hold true.

Testing \eqref{1.15} by $|b|^{p-2}b\ (2\leq p<\infty)$, using \eqref{1.12}, \eqref{1.28}, \eqref{1.29} and \eqref{1.31}, and setting $\phi=|b|^\frac{p}{2}$, and using the Gagliardo-Nirenberg inequality \cite{3}:
\begin{equation}
\|\phi\|_{L^{\frac{2s}{s-2},2}}\leq C\|\phi\|_{L^2}^{1-\frac3s}\|\phi\|_{H^1}^\frac3s\ \ \mathrm{with}\ \ 3<s\leq\infty,\label{4.5}
\end{equation}
and the generalized H\"{o}lder inequality \cite{13}:
\begin{equation}
\|fg\|_{L^{p,q}}\leq C\|f\|_{L^{p_1,q_1}}\|g\|_{L^{p_2,q_2}}\label{4.6}
\end{equation}
with $\frac1p=\frac{1}{p_1}+\frac{1}{p_2}$ and $\frac1q=\frac{1}{q_1}+\frac{1}{q_2}$, we derive
\begin{eqnarray*}
&&\frac1p\frac{d}{dt}\int_\Omega|b|^p dx+\frac12\int_\Omega|b|^{p-2}|\nabla b|^2dx+4\frac{p-2}{p^2}\int_\Omega|\nabla|b|^\frac{p}{2}|^2dx\\
&=&-\int_{\partial\Omega}|b|^{p-2}(b\cdot\nabla)\nu\cdot b d\sigma+\int_\Omega(b\cdot\nabla)u\cdot|b|^{p-2}b dx\\
&\leq&\|\nabla\nu\|_{L^\infty}\int_{\partial\Omega}|b|^p d\sigma-\sum\limits_i\int_\Omega b_iu\partial_i(|b|^{p-2}b)dx\\
&\leq&C\int_{\partial\Omega}\phi^2d\sigma+C\int_\Omega|u\phi\nabla\phi|dx\\
&\leq&C\int_{\partial\Omega}\phi^2d\sigma+C\|u\|_{L_w^s}\|\phi\|_{L^{\frac{2s}{s-2},2}}\|\nabla\phi\|_{L^2}\\
&\leq&C\|\phi\|_{L^2}\|\phi\|_{H^1}+C\|u\|_{L_w^s}\|\phi\|_{L^2}^{1-\frac3s}\|\nabla\phi\|_{L^2}^{1+\frac3s}\\
&\leq&2\frac{p-2}{p^2}\int_\Omega|\nabla\phi|^2dx+C\|\phi\|_{L^2}^2+C\|u\|_{L_w^s}^\frac{2s}{s-3}\|\phi\|_{L^2}^2,
\end{eqnarray*}
which yields
\begin{eqnarray*}
&&\frac{d}{dt}\int_\Omega\phi^2dx+C\int_\Omega|\nabla\phi|^2dx\leq C(1+\|u\|_{L_w^s}^\frac{2s}{s-3})\|\phi\|_{L^2}^2\\
&\leq&C\left(1+\frac{\|u\|_{L_w^s}^\frac{2s}{s-3}}{1+\log(e+\|u\|_{L_w^s})}\right)\|\phi\|_{L^2}^2(1+\log(e+\|u\|_{L_w^s}))\\
&\leq&C\left(1+\frac{\|u\|_{L_w^s}^\frac{2s}{s-3}}{1+\log(e+\|u\|_{L_w^s})}\right)(1+\log(e+y))\|\phi\|_{L^2}^2,
\end{eqnarray*}
from which it follows that
\begin{equation}
\|b\|_{L^\infty(t_0,t;L^p)}+\int_{t_0}^t\int_\Omega|b|^2|\nabla b|^2dx d\tau\leq C(e+y(t))^{C_0\epsilon}\label{4.7}
\end{equation}
with $$y(t):=\sup\limits_{[t_0,t]}\|u\|_{W^{1,4}}$$ for any $0<t_0\leq t\leq T$ and $C_0$ is an absolute constant, provided that
\begin{equation}
\int_{t_0}^T\frac{\|u\|_{L_w^s}^\frac{2s}{s-3}}{1+\log(e+\|u\|_{L_w^s})}d\tau\leq\epsilon<<1.\label{4.8}
\end{equation}

Testing \eqref{1.14} by $u_t$, using \eqref{1.12} and \eqref{1.13}, we infer that
\begin{eqnarray}
\frac12\frac{d}{dt}\int_\Omega|\nabla u|^2dx+\int_\Omega\rho|u_t|^2dx&=&-\int_\Omega\rho u\cdot\nabla u\cdot u_t dx+\int_\Omega b\cdot\nabla b\cdot u_t dx\nonumber\\
&=&:I_1+I_2.\label{4.9}
\end{eqnarray}

We first compute $I_2$.
\begin{eqnarray}
I_2&=&\int_\Omega\dv(b\otimes b)\cdot u_t dx=-\int_\Omega b\otimes b:\nabla u_t dx\nonumber\\
&=&-\frac{d}{dt}\int_\Omega b\otimes b:\nabla u dx+2\int_\Omega b\otimes b_t:\nabla u dx\nonumber\\
&\leq&-\frac{d}{dt}\int_\Omega b\otimes b:\nabla u dx+C\|b_t\|_{L^2}\|b\|_{L^6}\|\nabla u\|_{L^3}\nonumber\\
&\leq&-\frac{d}{dt}\int_\Omega b\otimes b:\nabla u dx+C\|b_t\|_{L^2}\|b\|_{L^6}\|\nabla u\|_{L^2}^\frac12\|u\|_{H^2}^\frac12\nonumber\\
&\leq&-\frac{d}{dt}\int_\Omega b\otimes b:\nabla u dx+\delta\|b_t\|_{L^2}^2+\delta\|u\|_{H^2}^2+C\|b\|_{L^6}^4\|\nabla u\|_{L^2}^2\label{4.10}
\end{eqnarray}
for any $0<\delta<1$.

We use \eqref{4.1}, \eqref{4.5} and \eqref{4.6} to bound $I_1$ as follows.
\begin{eqnarray}
I_1&\leq&\|\sqrt\rho u_t\|_{L^2}\|\sqrt\rho\|_{L^\infty}\|u\|_{L_w^s}\|\nabla u\|_{L^{\frac{2s}{s-2},2}}\nonumber\\
&\leq&C\|\sqrt\rho u_t\|_{L^2}\|u\|_{L_w^s}\|\nabla u\|_{L^2}^{1-\frac3s}\|u\|_{H^2}^\frac3s\nonumber\\
&\leq&\delta\|\sqrt u_t\|_{L^2}^2+\delta\|u\|_{H^2}^2+C\|u\|_{L_w^s}^\frac{2s}{s-3}\|\nabla u\|_{L^2}^2\label{4.11}
\end{eqnarray}
for any $0<\delta<1$.

On the other hand, by the $H^2$-theory of the Stokes system, and using \eqref{4.1}, \eqref{4.5} and \eqref{4.6}, we obtain
\begin{eqnarray*}
\|u\|_{H^2}&\leq&C\left\|-\Delta u+\nabla\left(\pi+\frac12|b|^2\right)\right\|_{L^2}\\
&\leq&C\|\rho\partial_t u+\rho u\cdot\nabla u-b\cdot\nabla b\|_{L^2}\\
&\leq&C\|\sqrt\rho u_t\|_{L^2}+C\|u\|_{L_w^s}\|\nabla u\|_{L^{\frac{2s}{s-2},2}}+C\|b\cdot\nabla b\|_{L^2}\\
&\leq&C\|\sqrt\rho u_t\|_{L^2}+C\|u\|_{L_w^s}\|\nabla u\|_{L^2}^{1-\frac3s}\|u\|_{H^2}^\frac3s+C\|b\cdot\nabla b\|_{L^2},
\end{eqnarray*}
which gives
\begin{equation}
\|u\|_{H^2}\leq C\|\sqrt\rho u_t\|_{L^2}+C\|b\cdot\nabla b\|_{L^2}+C\|u\|_{L_w^s}^\frac{s}{s-3}\|\nabla u\|_{L^2}.\label{4.12}
\end{equation}

Testing \eqref{1.15} by $b_t-\Delta b$, using \eqref{4.5} and \eqref{4.6}, we deduce that
\begin{eqnarray}
&&\frac{d}{dt}\int_\Omega|\cl b|^2dx+\int_\Omega(|b_t|^2+|\Delta b|^2)dx\nonumber\\
&=&\int_\Omega(b\cdot\nabla u-u\cdot\nabla b)(b_t-\Delta b)dx\nonumber\\
&\leq&(\|u\|_{L_w^s}\|\nabla b\|_{L^{\frac{2s}{s-2},2}}+\|b\|_{L^6}\|\nabla u\|_{L^3})(\|b_t\|_{L^2}+\|\Delta b\|_{L^2})\nonumber\\
&\leq&C(\|u\|_{L_w^s}\|\nabla b\|_{L^2}^{1-\frac3s}\|b\|_{H^2}^\frac3s+C\|b\|_{L^6}\|\nabla u\|_{L^2}^\frac12\|u\|_{H^2}^\frac12)(\|b_t\|_{L^2}+\|\Delta b\|_{L^2})\nonumber\\
&\leq&\frac12(\|b_t\|_{L^2}^2+\|\Delta b\|_{L^2}^2)+\delta\|u\|_{H^2}^2+C\|b\|_{L^6}^4\|\nabla u\|_{L^2}^2+C\|u\|_{L_w^s}^\frac{2s}{s-3}\|\nabla b\|_{L^2}^2+C\label{4.13}
\end{eqnarray}
for any $0<\delta<1$.

It is easy to compute that
\begin{eqnarray}
\frac{d}{dt}\int_\Omega|b|^4dx&\leq&C\int_\Omega|b|^3|b_t|dx\nonumber\\
&\leq&C\|b\|_{L^6}^3\|b_t\|_{L^2}\leq\delta\|b_t\|_{L^2}^2+C\|b\|_{L^6}^6\label{4.14}
\end{eqnarray}
for any $0<\delta<1$.

Combining \eqref{4.9}, \eqref{4.10}, \eqref{4.11}, \eqref{4.12}, \eqref{4.13} and \eqref{4.14}, and taking $\delta$ small enough, we obtain
\begin{eqnarray}
&&\frac{d}{dt}\int_\Omega(|\nabla u|^2+|\cl b|^2+b\otimes b:\nabla u+C_0|b|^4)dx+\int_\Omega(\rho|u_t|^2+|b_t|^2+|\Delta b|^2)dx+\|u\|_{H^2}^2\nonumber\\
&\leq&C\|b\|_{L^6}^4\|\nabla u\|_{L^2}^2+C\|u\|_{L_w^s}^\frac{2s}{s-3}(\|\nabla u\|_{L^2}^2+\|\cl b\|_{L^2}^2)+C\|b\cdot\nabla b\|_{L^2}^2+C.\label{4.15}
\end{eqnarray}

Using \eqref{4.4}, \eqref{4.7}, \eqref{4.8} and the Gronwall inequality, we have
\begin{eqnarray}
&&\int_\Omega(|\nabla u|^2+|\cl b|^2+b\otimes b:\nabla u+C_0|b|^4)dx\nonumber\\
&\leq&\left[\int_\Omega(|\nabla u_0|^2+|\cl b_0|^2+b_0\otimes b_0:\nabla u_0+C_0|b_0|^4)dx\right.\nonumber\\
&&\left.+C\|b\|_{L^\infty(t_0,t;L^6)}^4\int_{t_0}^t\|\nabla u\|_{L^2}^2d\tau+C(t-t_0)+C\int_{t_0}^t\|b\cdot\nabla b\|_{L^2}^2d\tau\right]\exp\left(\int_{t_0}^t\|u\|_{L_w^s}^\frac{2s}{s-3}d\tau\right)\nonumber\\
&\leq&C(e+y)^{C_0\epsilon}\exp\left[\int_{t_0}^t\frac{\|u\|_{L_w^s}^\frac{2s}{s-3}}{1+\log(e+\|u\|_{L_w^s})} d\tau(1+\log(e+y))\right]\nonumber\\
&\leq&C(e+y)^{C_0\epsilon}.\label{4.16}
\end{eqnarray}

Plugging \eqref{4.16} into \eqref{4.15} and integrating over $[t_0,t]$, we have
\begin{equation}
\int_{t_0}^t\int_\Omega(\rho|u_t|^2+|b_t|^2+|\Delta b|^2)dx d\tau+\int_{t_0}^t\|u\|_{H^2}^2d\tau\leq C(e+y)^{C_0\epsilon}.\label{4.17}
\end{equation}

Applying $\partial_t$ to \eqref{1.15}, testing by $u_t$, using \eqref{1.12} and \eqref{1.13}, we obtain
\begin{eqnarray}
&&\frac12\frac{d}{dt}\int_\Omega\rho|u_t|^2dx+\int_\Omega|\nabla u_t|^2dx\nonumber\\
&=&-\int_\Omega\rho u\cdot\nabla|u_t|^2dx-\int_\Omega\rho u\cdot\nabla(u\cdot\nabla u\cdot u_t)dx\nonumber\\
&&-\int_\Omega\rho u_t\cdot\nabla u\cdot u_t dx+\int_\Omega b\otimes b_t:\nabla u_t dx+\int_\Omega b_t\otimes b:\nabla u_t dx\nonumber\\
&\leq&C\|u\|_{L^6}\|\sqrt\rho u_t\|_{L^3}\|\nabla u_t\|_{L^2}+C\|u\|_{L^6}\|\nabla u\|_{L^6}\|u_t\|_{L^6}\|\nabla u\|_{L^2}\nonumber\\
&&+C\|u\|_{L^6}^2\|\Delta u\|_{L^2}\|u_t\|_{L^6}+C\|u\|_{L^6}^2\|\nabla u\|_{L^6}\|\nabla u_t\|_{L^2}\nonumber\\
&&+C\|\sqrt\rho u_t\|_{L^4}^2\|\nabla u\|_{L^2}+C\|b\|_{L^6}\|b_t\|_{L^3}\|\nabla u_t\|_{L^2}\nonumber\\
&\leq&C\|\nabla u\|_{L^2}\|\sqrt\rho u_t\|_{L^2}^\frac12\|\sqrt\rho u_t\|_{L^6}^\frac12\|\nabla u_t\|_{L^2}\nonumber\\
&&+C\|\nabla u\|_{L^2}^2\|u\|_{H^2}\|\nabla u_t\|_{L^2}+C\|\nabla u\|_{L^2}\|\sqrt\rho u_t\|_{L^2}^\frac12\|\sqrt\rho u_t\|_{L^6}^\frac32\nonumber\\
&&+C\|b\|_{L^6}\|b_t\|_{L^3}\|\nabla u_t\|_{L^2}\nonumber\\
&\leq&C\|\nabla u\|_{L^2}\|\sqrt\rho u_t\|_{L^2}^\frac12\|\nabla u_t\|_{L^2}^\frac32+C\|\nabla u\|_{L^2}^2\|u\|_{H^2}\|\nabla u_t\|_{L^2}\nonumber\\
&&+C\|\nabla u\|_{L^2}\|\sqrt\rho u_t\|_{L^2}^\frac12\|\nabla u_t\|_{L^2}^\frac32+C\|b\|_{L^6}\|b_t\|_{L^3}\|\nabla u_t\|_{L^2}\nonumber\\
&\leq&\frac14\|\nabla u_t\|_{L^2}^2+C\|\nabla u\|_{L^2}^4(\|\sqrt\rho u_t\|_{L^2}^2+\|u\|_{H^2}^2)+C\|b\|_{L^6}^2\|b_t\|_{L^3}^2\nonumber\\
&\leq&\frac14\|\nabla u_t\|_{L^2}^2+C\|\nabla u\|_{L^2}^4(\|\sqrt\rho u_t\|_{L^2}^2+\|u\|_{H^2}^2)+\frac14\|\cl b_t\|_{L^2}^2+C\|b\|_{L^6}^4\|b_t\|_{L^2}^2.\label{4.18}
\end{eqnarray}

Applying $\partial_t$ to \eqref{1.15}, testing by $b_t$ and using \eqref{1.12}, we get
\begin{eqnarray}
&&\frac12\frac{d}{dt}\int_\Omega|b_t|^2dx+\int_\Omega|\cl b_t|^2dx\nonumber\\
&=&-\int_\Omega(u_t\cdot\nabla b-b_t\nabla u-b\cdot\nabla u_t)b_t dx\nonumber\\
&\leq&\|u_t\|_{L^6}\|\nabla b\|_{L^2}\|b_t\|_{L^3}+\|\nabla u\|_{L^2}\|b_t\|_{L^4}^2+\|\nabla u_t\|_{L^2}\|b\|_{L^6}\|b_t\|_{L^3}\nonumber\\
&\leq&\frac14\|\nabla u_t\|_{L^2}^2+\frac14\|\cl b_t\|_{L^2}^2+C\|\nabla b\|_{L^2}^4\|b_t\|_{L^2}^2+C\|\nabla u\|_{L^2}^4\|b_t\|_{L^2}^2.\label{4.19}
\end{eqnarray}

Combining \eqref{4.18} and \eqref{4.19} and integrating over $[t_0,t]$, we have
\begin{equation}
\int_\Omega(|\rho u_t|^2+|b_t|^2)dx+\int_{t_0}^t\int_\Omega(|\nabla u_t|^2+|\cl b_t|^2)dx d\tau\leq C(e+y)^{C_0\epsilon}.\label{4.20}
\end{equation}

Similarly to \eqref{4.12}, we deduce that
\begin{eqnarray*}
\|u\|_{H^2}&\leq&C\|\sqrt\rho u_t\|_{L^2}+C\|u\|_{L^6}\|\nabla u\|_{L^3}+C\|b\|_{L^6}\|\nabla b\|_{L^3}\\
&\leq&C\|\sqrt\rho u_t\|_{L^2}+C\|u\|_{L^6}\|\nabla u\|_{L^2}^\frac12\|u\|_{H^2}^\frac12+C\|b\|_{L^6}\|\nabla b\|_{L^2}^\frac12\|b\|_{H^2}^\frac12,
\end{eqnarray*}
which leads to
\begin{equation}
\|u\|_{H^2}^2\leq C\|\sqrt\rho u_t\|_{L^2}^2+C\|\nabla u\|_{L^2}^6+C\|\nabla b\|_{L^2}^6+\frac12\|b\|_{H^2}^2.\label{4.21}
\end{equation}

Similarly, we have
\begin{eqnarray*}
\|b\|_{H^2}&\leq&C\|b_t+u\cdot\nabla b-b\cdot\nabla u\|_{L^2}\\
&\leq&C\|b_t\|_{L^2}+C\|u\|_{L^6}\|\nabla b\|_{L^3}+C\|b\|_{L^6}\|\nabla u\|_{L^3}\\
&\leq&C\|b_t\|_{L^2}+C\|u\|_{L^6}\|\nabla b\|_{L^2}^\frac12\|b\|_{H^2}^\frac12+C\|b\|_{L^6}\|\nabla u\|_{L^2}^\frac12\|u\|_{H^2}^\frac12,
\end{eqnarray*}
which implies
\begin{equation}
\|b\|_{H^2}^2\leq C\|b_t\|_{L^2}^2+C\|\nabla u\|_{L^2}^6+C\|\nabla b\|_{L^2}^6+\frac12\|u\|_{H^2}^2.\label{4.22}
\end{equation}

Combining \eqref{4.21} and \eqref{4.22}, using \eqref{4.20} and \eqref{4.16}, we conclude that
\begin{equation}
\|u\|_{H^2}^2+\|b\|_{H^2}^2\leq C(e+y)^{C_0\epsilon},\label{4.23}
\end{equation}
and thus
\begin{equation}
\|u\|_{L^\infty(0,T;H^2)}+\|b\|_{L^\infty(0,T;H^2)}\leq C.\label{4.24}
\end{equation}

Now it is standard to prove that
\begin{eqnarray}
&&\|u\|_{L^2(0,T;H^3)}+\|b\|_{L^2(0,T;H^3)}\leq C,\label{4.25}\\
&&\|\rho\|_{L^\infty(0,T;H^2)}\leq C.\label{4.26}
\end{eqnarray}

\noindent (II) Let \eqref{1.22} hold true.

Similarly to \eqref{4.7}, we take $s=\infty$ and using \eqref{1.30}, we still get \eqref{4.7} provided that
\begin{equation}
\int_{t_0}^T\|u(t)\|_{BMO}^2dt\leq\epsilon<<1.\label{4.27}
\end{equation}

We still have \eqref{4.9}, \eqref{4.10}, \eqref{4.11} with $s=\infty$, \eqref{4.12} with $s=\infty$, \eqref{4.13} with $s=\infty$, \eqref{4.14}, \eqref{4.15} with $s=\infty$, and then using \eqref{4.27} and \eqref{1.30}, we arrive at \eqref{4.16} and \eqref{4.17}. Then by the same calculations as that in \eqref{4.18}-\eqref{4.26}, we conclude that \eqref{4.18}-\eqref{4.26} hold true.

This completes the proof of Theorem \ref{th1.3}.

\hfill$\square$

%%%%%%%%%%%%%%%----------------------------------------------------------------------------------------------
\section{Proof of Theorem \ref{th1.4}}

This section is devoted to the proof of Theorem \ref{th1.4}, we only need to establish a priori estimates.

First, using the formula $a\times(b\times c)=(a\cdot c)b-(a\cdot b)c$, and the fact that $|d|=1$ implies $d\Delta d=-|\nabla d|^2$, we have the following equivalent equation
\begin{equation}
\frac12 d_t-\frac12 d\times d_t=\Delta d+d|\nabla d|^2.\label{5.1}
\end{equation}

Testing \eqref{5.1} by $d_t$ and using $(a\times b)\cdot b=0$ and $d\cdot d_t=0$, we get
\begin{equation}
\frac{d}{dt}\int_\Omega|\nabla d|^2dx+\int_\Omega|d_t|^2dx\leq 0.\label{5.2}
\end{equation}

Testing \eqref{1.23} by $-\Delta d_t$ and using $|d|=1$, we find that
\begin{eqnarray}
&&\frac12\frac{d}{dt}\int_\Omega|\Delta d|^2dx+\int_\Omega|\nabla d_t|^2dx=-\int_\Omega(d|\nabla d|^2+d\times\Delta d)\cdot\Delta d_t dx\nonumber\\
&=&\int_\Omega\nabla(d|\nabla d|^2+d\times\Delta d)\cdot\nabla d_t dx\nonumber\\
&\leq&C(\|\nabla d\|_{L^q}\|\nabla d\|_{L^\frac{4q}{q-2}}^2+\|\nabla d\|_{L^q}\|\Delta d\|_{L^\frac{2q}{q-2}}+\|\nabla\Delta d\|_{L^2})\|\nabla d_t\|_{L^2}\nonumber\\
&\leq&C(\|\nabla d\|_{L^q}\|\Delta d\|_{L^\frac{2q}{q-2}}+\|\nabla\Delta d\|_{L^2})\|\nabla d_t\|_{L^2}\nonumber\\
&\leq&C(\|\nabla d\|_{L^q}\|\Delta d\|_{L^2}^{1-\frac3q}\|d\|_{H^3}^\frac3q+\|d\|_{H^3})\|\nabla d_t\|_{L^2}\nonumber\\
&\leq&\frac14\|\nabla d_t\|_{L^2}^2+\delta\|d\|_{H^3}^2+C\|\nabla d\|_{L^q}^\frac{2q}{q-3}\|\Delta d\|_{L^2}^2\label{5.3}
\end{eqnarray}
for any $0<\delta<1$. Here we have used the Gagliardo-Nirenberg inequalities:
\begin{eqnarray}
&&\|\nabla d\|_{L^\frac{4q}{q-2}}^2\leq C\|d\|_{L^\infty}\|\Delta d\|_{L^\frac{2q}{q-2}},\label{5.4}\\
&&\|\Delta d\|_{L^\frac{2q}{q-2}}\leq C\|\Delta d\|_{L^2}^{1-\frac3q}\|d\|_{H^3}^\frac3q.\label{5.5}
\end{eqnarray}

Applying $\partial_i$ to \eqref{1.23}, we get $$\partial_id_t-\Delta\partial_id=\partial_i(d|\nabla d|^2)+\partial_id\times\Delta d+d\times\Delta\partial_id.$$

Testing the above equation by $\Delta\partial_id$, summing over $i$, and using \eqref{5.4} and \eqref{5.5} and $|d|=1$, we obtain
\begin{eqnarray*}
\|d\|_{H^3}&\leq&C(\|d\|_{L^2}+\|\nabla\Delta d\|_{L^2})\\
&\leq&C+C\|\nabla d_t\|_{L^2}+C\|\nabla(d|\nabla d|^2)\|_{L^2}+\sum\limits_iC\|\partial_id\times\Delta d\|_{L^2}\\
&\leq&C+C\|\nabla d_t\|_{L^2}+C\|\nabla d\|_{L^q}\|\nabla d\|_{L^\frac{4q}{q-2}}^2+C\|\nabla d\|_{L^q}\|\Delta d\|_{L^\frac{2q}{q-2}}\\
&\leq&C+C\|\nabla d_t\|_{L^2}+C\|\nabla d\|_{L^q}\|\Delta d\|_{L^\frac{2q}{q-2}}\\
&\leq&C+C\|\nabla d_t\|_{L^2}+C\|\nabla d\|_{L^q}\|\Delta d\|_{L^2}^{1-\frac3q}\|d\|_{H^3}^\frac3q,
\end{eqnarray*}
which yields
\begin{equation}
\|d\|_{H^3}\leq C+C\|\nabla d_t\|_{L^2}+C\|\nabla d\|_{L^q}^\frac{q}{q-3}\|\Delta d\|_{L^2}.\label{5.6}
\end{equation}

Plugging \eqref{5.6} into \eqref{5.3} and taking $\delta$ small enough, we have
\begin{eqnarray*}
&&\frac{d}{dt}\int_\Omega|\Delta d|^2dx+\int_\Omega|\nabla d_t|^2dx\\
&\leq&C+C\|\nabla d\|_{L^q}^\frac{2q}{q-3}\|\Delta d\|_{L^2}^2\\
&\leq&C+C\frac{\|\nabla d\|_{L^q}^\frac{2q}{q-3}}{1+\log(e+\|\nabla d\|_{L^q})}\|\Delta d\|_{L^2}^2\log(e+\|\nabla d\|_{L^q})\\
&\leq&C+C\frac{\|\nabla d\|_{L^q}^\frac{2q}{q-3}}{1+\log(e+\|\nabla d\|_{L^q})}\|\Delta d\|_{L^2}^2\log(e+y),
\end{eqnarray*}
which implies
\begin{equation}
\int_\Omega|\Delta d|^2dx+\int_{t_0}^t\int_\Omega|\nabla d_t|^2dx d\tau\leq C(e+y)^{C_0\epsilon}.\label{5.7}
\end{equation}

Provided that $$\int_{t_0}^T\frac{\|\nabla d\|_{L^q}^\frac{2q}{q-3}}{1+\log(e+\|\nabla d\|_{L^q})}d\tau\leq\epsilon<<1,$$ with $y(t):=\sup\limits_{[t_0,t]}\|d\|_{H^3}$ for any $0<t_0\leq t\leq T$ and $C_0$ is an absolute constant.

It follows from \eqref{1.23}, \eqref{5.6} and \eqref{5.7} that
\begin{equation}
\int_\Omega|d_t|^2dx+\int_{t_0}^t\|d\|_{H^3}^2d\tau\leq C(e+y)^{C_0\epsilon}.\label{5.8}
\end{equation}

Applying $\partial_t$ to \eqref{1.23}, testing by $-\Delta d_t$, and using $|d|=1$, \eqref{5.7} and \eqref{5.8}, we have
\begin{eqnarray*}
&&\frac12\frac{d}{dt}\int_\Omega|\nabla d_t|^2dx+\int_\Omega|\Delta d_t|^2dx=-\int_\Omega[\partial_t(d|\nabla d|^2)+d_t\times\Delta d]\Delta d_t dx\\
&\leq&C(\|\nabla d\|_{L^6}^2\|d_t\|_{L^6}+\|\nabla d\|_{L^6}\|\nabla d_t\|_{L^3}+\|d_t\|_{L^\infty}\|\Delta d\|_{L^2})\|\Delta d_t\|_{L^2}\\
&\leq&C(\|\nabla d\|_{L^6}^2\|d_t\|_{L^6}+\|\Delta d\|_{L^2}\|\nabla d_t\|_{L^2}^\frac12\|\Delta d_t\|_{L^2}^\frac12+\|\Delta d\|_{L^2}\|d_t\|_{L^2})\|\Delta d_t\|_{L^2}\\
&\leq&\frac12\|\Delta d_t\|_{L^2}^2+C\|d\|_{H^2}^4\|d_t\|_{H^1}^2+C\|d\|_{H^2}^2\|d_t\|_{L^2}^2,
\end{eqnarray*}
which implies
\begin{equation}
\int_\Omega|\nabla d_t|^2dx+\int_{t_0}^t\|\Delta d_t\|_{L^2}^2d\tau\leq C(e+y)^{C_0\epsilon}.\label{5.9}
\end{equation}

It follows from \eqref{5.6}, \eqref{5.7}, \eqref{5.8} and \eqref{5.9} that $$\|d\|_{H^3}\leq C+C\|\nabla d_t\|_{L^2}+C\|\nabla d\|_{L^6}^2\|\Delta d\|_{L^2}\leq C(e+y)^{C_0\epsilon},$$ which leads to $$\|d\|_{L^\infty(0,T;H^3)}\leq C.$$

This completes the proof of Theorem \ref{th1.4}.

\hfill$\square$
%%%%%%%%%%%%%%%----------------------------------------------------------------------------------------------
\section{Acknowledgments}

J. Fan is partially supported by NSFC (No. 11171154), Junpin Yin is supported by the NSFC (Grant No.111 01044) and Beijing Center for Mathematics and Information
Interdisciplinary Sciences (BCMIIS). The authors would like to thank the referee
for reading the paper carefully and for the valuable comments which improved the presentation of the paper.

%%%%%%%%%%%%%%----------------------------------------------------------------------------------------------

\end{document}